\documentclass[10pt,a4paper,DIV12]{scrartcl}         

\usepackage[utf8]{inputenc}
\usepackage[T1]{fontenc}
\usepackage{lmodern}
\usepackage{amsmath,amsfonts, amssymb, mathrsfs}
\usepackage{bbm}
\usepackage{graphicx}
\usepackage{trfsigns}
\usepackage{algpseudocode} 
\usepackage{algorithm}
\usepackage{booktabs}
\usepackage{listings}
\usepackage{xcolor}
\usepackage{overpic}
\usepackage{url}
\usepackage{stmaryrd}
\usepackage{subcaption}
\numberwithin{equation}{section}

\begin{document}
\title{Computational comparison of surface metrics for PDE constrained shape
optimization}
\author{Volker Schulz\thanks{Universit\"at Trier, Universit\"atsring 15, D-54296 Trier, Germany, Email: volker.schulz@uni-trier.de, siebenborn@uni-trier.de} \and Martin Siebenborn\footnotemark[1]
}
\date{}
\maketitle

\newcommand{\Gv}{{\Gamma}}
\newcommand{\Gi}{{\Gamma_\text{in}}}
\newcommand{\Go}{{\Gamma_\text{out}}}
\newcommand{\Gw}{{\Gamma_\text{wall}}}
\newcommand{\Oe}{{\Omega_\text{ext}}}
\newcommand{\N}{{\mathbbm{N}}} 
\newcommand{\R}{{\mathbbm{R}}} 
\newcommand{\SM}{{\cal B}_i^{1/2}(\Gamma_0)} 

\def\scp#1{\left\langle #1\right\rangle}
\def\norm#1{\left\| #1\right\|}
\def\snorm#1{\| #1\|}
\def\CROP#1{}
\def\grad{\mbox{grad}}
\newcommand{\LL}{{\mathscr{L}}}
\newcommand{\HH}{{\mathscr{H}}}

\def\bei#1{\vrule width 0.4pt height 14pt depth 9pt
           \lower 8pt \hbox{$ _{\hbox{} #1}$}\!\!\!}

\newcommand{\ACS}{2pt}
\newcommand{\mat}[4]{{\arraycolsep\ACS
\left#1\begin{array}{@{}*{#2}{c}@{}}#4\end{array}\right#3}}

\begin{abstract}
We compare surface metrics for shape optimization problems with constraints, consisting mainly of partial differential equations (PDE), from a computational point of view.
In particular, classical Laplace-Beltrami type based metrics are compared with Steklov-Poincaré type metrics.
The test problem is the minimization of energy dissipation of a body in a Stokes flow.
We therefore set up a quasi-Newton method on appropriate shape manifolds together with an augmented Lagrangian framework, in order to enable a straightforward integration of geometric constraints for the shape.
The comparison is focussed towards convergence behavior as well as
effects on the mesh quality during shape optimization.
\end{abstract}


\section{Introduction}
Shape optimization is a challenging field with many interesting applications. As examples, we mention aerodynamic shape optimization \cite{AIAA-2013}, acoustic shape optimization \cite{Berggren-horn-2007} or optimization of interfaces in transmission problems \cite{Langer-2015,Skin-2015,Paganini}. The general structure of such a a PDE constrained shape optimization problem is of the form
\begin{align*}
\min\limits_{y,\Omega} J(y,&\Omega)\\
\mbox{s.t.}\quad \mbox{B}_\Omega(y)&=0\\
 c(y,\Omega) &= 0
\end{align*}
where $\Omega\subset D$ is an open subset of a hold-all domain $D\subset\R^d$ and $J$ is a real valued functional. Usually only some part of the boundary $\Gamma\subset\partial\Omega$ is free for optimization. The exact description of the set from which $\Gamma$ is taken, is given in section 3. The constraint $\mbox{B}_\Omega(y)=0$ denotes a boundary value problem defined on the domain $\Omega$ given in the form of equations in appropriate function spaces, where $y$ is the solution of the boundary value problem $B_\Omega$ consisting of a set of partial differential equations together with some boundary conditions. Furthermore $c(y,\Omega)$ denotes a finite number of sufficiently smooth constraints. The additional constraints maybe even in the form of inequalities, but this aspect is not in the focus of this paper.  

In principle, there are two major conceptual approaches: the direct parametrization approach, which a priorily parameterizes the shape $\Gamma$ to be optimized, e.g., within a CAD framework, and the shape calculus approach operating in shape spaces. The direct parameterization approach suffers from obvious limitations with respect to the reachable geometries, but can be embedded within a vector space framework which simplifies the numerical treatment significantly and thus enables the application of classical methods for PDE constrained optimization \cite{optbook}. The shape calculus approach has received strong attention in particular concerning its theoretical framework in the form of the calculus for generating shape derivative information. Several books and publications deal with this aspect in detail, e.g. \cite{Delfour-Zolesio-2001,SokoZol}. In contrast to that, only relatively few publications study computational aspects of shape optimization based on the shape calculus. In particular, studies on proper shape metrics for usage within gradient-type methods and also on the efficient numerical treatment of additional constraints are missing. In most cases, steepest descent methods are applied involving representations of the shape derivative, which are based on the $L^2$ scalar product or a Laplace-Beltrami scalar product or combinations thereof. In \cite{SIOPT2015}, a novel scalar product based on the Steklov-Poincar\'e operator is introduced in the context of shape optimization and shown to posses the following advantadgeous properties: volumetric and boundary expressions of the shape derivative can be treated in a consistent manner; and the resulting shape manifold admits kinks and is complete. Furthermore, this scalar product is in line with the coercivity results for shape Hessians in \cite{EHS2007} for elliptic problems.

Therefore, we compare in this paper the novel metric introduced in \cite{SIOPT2015} with the metric used so far in many publications from a computational perspective. The test case is the computation of a shape embedded in a Stokes-flow which minimizes drag and satisfies further geometric constraints. The optimal solution is well-known as the so called Haack ogive \cite{Haack1941,Pironneau1973}. The geometric constraints are taken into account within an augmented Lagrangian framework similarly to \cite{EHS2007}. Furthermore, quasi-Newton techniques are applied in order to accelerate the converge of the shape optimization scheme. A major issue of the comparison of metrics is the surface and volume mesh quality of the deformed computational mesh. We observe that the mesh quality only mildly deteriorates during the optimization iterations based on the Steklov-Poincar\'e metric, although the overall deformation is considerably large. In contrast to that, the mesh quality drastically deteriorates during iterations with standard metrics.

This paper is organized in the following way. Section 2 gives details on the specific formulation used as a benchmark for the comparison of metrics. In section 3, we introduce the manifold point of view on shape optimization as well as the two specific metrics to be compared in later sections. Section 4 introduces an augmented Lagrangian formulation for an efficient treatment of the geometric constraints involved in the test case. Finally, the results on the computational comparison of the two types of metrics for shape optimization are presented in section 5. Conclusions are drawn in section 6. 

\section{Problem formulation: Optimal shapes in Stokes flows}
\begin{figure}
\begin{center}
\def\svgwidth{0.6\textwidth}
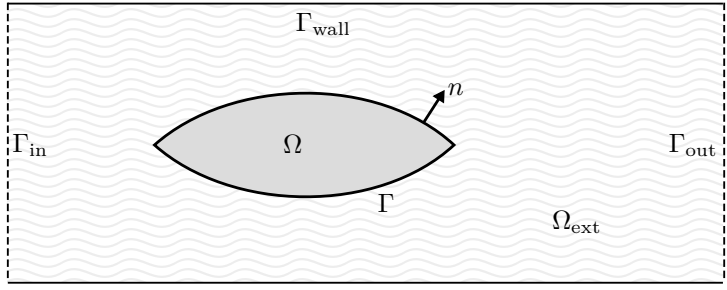
\end{center}
\caption{Schematic view of the flow field and the variable shape $\Omega$ encircled by $\Gv$}
\label{fig_schematic}
\end{figure}

We consider incompressible flow which is dominated by viscous forces around an obstacle described by the Stokes equations.
The aim is to shape a $d$-dimensional body such that the energy dissipation of the system is minimized under certain geometrical constraints.
For the optimization to be reasonable, the volume and the barycenter of the body are required to be constant.
Here, the dimension of the problem is fixed to $d=2$ or $d=3$.
This situation is visualized in Figure \ref{fig_schematic}, where $\Omega \subset \mathbbm{R}^d$ is the obstacle and $\Gv$ its boundary, which is considered to be variable.
$\Oe \subset \mathbbm{R}^d$ denotes the flow field.
This is the domain for the finite elements, whereas $\Omega$ is a hole in the discretization mesh.
For the geometric restrictions we need to compute the body's volume
\begin{equation}
\text{vol}(\Omega) = \int_\Omega 1\, dx\, \in \mathbbm{R}
\end{equation}
and its barycenter
\begin{equation}
\text{bc}^{\Omega} = \frac{1}{\text{vol}(\Omega)} \int_\Omega x\, dx\, \in \mathbbm{R}^d.
\end{equation}
Finally, we end up with the PDE constraint shape optimization problem
\begin{equation}\label{eq_objective}
\min\limits_{(v, \Omega)} J(v, \Omega) = \int_\Oe \sum\limits_{i,j=1}^{d} \left( \frac{\partial v_i}{\partial x_j} \right)^2\, dx
\end{equation}
subject to the Stokes equations, where the viscosity is normalized to 1
\begin{equation}\label{eq_stokes}
\begin{aligned}
\Delta v + \nabla p &= -f       &\quad& \text{in}\; \Oe\\
\text{div}\, v      &=  0       &     & \text{in}\; \Oe\\
                  v &= v_\infty &     & \text{on}\; \Gi \cup \Go\\
                  v &= 0        &     & \text{on}\; \Gw \cup \Gv\\
\end{aligned}
\end{equation}
and the geometric constraints, i.e.\ barycenter and volume
\begin{equation}\label{eq_geo_constraints}
\begin{aligned}
c_1(\Omega) &=\text{bc}(\Omega)_1 - \text{bc}(\Omega_0)_1\\
&\vdots\\
c_d(\Omega) &=\text{bc}(\Omega)_d - \text{bc}(\Omega_0)_d\\
c_{d+1}(\Omega) &=\text{vol}(\Omega) - \text{vol}(\Omega_0).
\end{aligned}
\end{equation}
In equation \eqref{eq_stokes} $v:\Oe \to \mathbbm{R}^d$ denotes the velocity and $p: \Oe \to \mathbbm{R}$.
For the corresponding weak formulation of the Stokes equation we assume $v \in H^1 (\Oe)^d$ and $p \in L_{2,0} (\Oe) := \lbrace q \in L_2 (\Oe) : \int_\Oe q\, dx = 0\rbrace$.
Since we neglect body forces like gravity we can set $f=0$.

For a gradient based optimization we need derivatives with respect to the shape which are defined in the following way.
The shape derivative in direction of a smooth vector field $V: \Omega \to \mathbbm{R}^d$ is defined as
\begin{equation}
d J (\Omega) \left[ V \right] := \lim\limits_{h \to 0+} \frac{J(\Omega_h) - J(\Omega)}{h}
\end{equation}
where $\Omega_h = \lbrace x + h\cdot V(x) \,:\, x \in \Omega \rbrace$ is perturbed according to $V$.
For the particular setting in \eqref{eq_objective} and \eqref{eq_stokes} the derivative of $J$ subject to Stokes equation is given by (cf.\ \cite{Mohammadi-2001})
\begin{equation}\label{obj-der}
dJ(\Omega)[V] = -\int_\Gv \langle n , V \rangle \sum\limits_{i=1}^d \left( \frac{\partial v_i}{\partial n} \right)^2\, ds.
\end{equation}
The derivatives of the geometric constraints can be derived by applying the calculus developed in \cite{Delfour-Zolesio-2001} to $c_1,\dots, c_{d+1}$ yielding
\begin{equation}
\begin{aligned}
dc_{i}(\Omega)[V] &= -\frac{1}{\left(\text{vol}(\Omega)\right)^2}   \int_\Gv \langle n, V \rangle \, ds   \int_\Omega x_i \, dx    +    \frac{1}{\text{vol}(\Omega)} \int_\Gv x_i \langle n, V \rangle \, ds\\
&= \frac{1}{\text{vol}(\Omega)} \int_\Gv \left( x_i -  \text{bc}(\Omega)_i\right) \langle n, V \rangle  \, ds
\end{aligned}
\end{equation}
and
\begin{equation}
dc_{d+1}(\Omega)[V] = \int_\Gv \langle n, V \rangle \, ds.
\end{equation}
In order to formulate an augmented Lagrangian method for shape optimization in section \ref{sec_augmented_lagrange}, derivatives of the squared constraints are also required.
By applying the chain rule we obtain
\begin{equation}
dc_{i}^2(\Omega)[V] = 2\left( \text{bc}(\Omega)_i - \text{bc}(\Omega_0)_i \right) \frac{1}{\text{vol}(\Omega)} \int_\Gv \left( x_i -  \text{bc}(\Omega)_i\right) \langle n, V \rangle  \, ds
\end{equation}
and
\begin{equation}
dc_{d+1}^2(\Omega)[V] = 2\left( \text{vol}(\Omega) - \text{vol}(\Omega_0) \right) \int_\Gv \langle n, V \rangle \, ds.
\end{equation}

\section{Metrics on the manifold of feasible shapes}

In \cite{VHS-shape-Riemann}, it is pointed out that shape optimization can be viewed as optimization on Riemannian shape manifolds and resulting optimization methods can be constructed and analyzed within this framework, which combines algorithmic ideas from \cite{Absil-book-2008} with the differential geometric point of view established in \cite{MM-2006}.
Let us study connected and compact curves $\Gamma$ as in figure \ref{fig_schematic}. Although the optimal solution of our test problem has two  kinks, we first consider smooth closed curves in order to discuss metrics more easily. 

In \cite{MM-2006}, this set of smooth closed curves is characterized by  
\begin{equation}
B_e(S^1,\R^2):=\mbox{Emb}(S^1,\R^2)/\mbox{Diff}(S^1),
\end{equation}
i.e., as the set of all equivalence classes of $C^\infty$ embeddings of $S^1$ into the plane ($\mbox{Emb}(S^1,\R^2)$), where the equivalence relation is defined by the set of all $C^\infty$ re-parameterizations, i.e., diffeomorphisms of $S^1$ into itself ($\mbox{Diff}(S^1)$). A particular point on the manifold $B_e(S^1,\R^2)$ is represented by a curve $\Gamma\colon S^1\ni\theta\mapsto \Gamma(\theta)\in\R^2$. Because of the equivalence relation ($\mbox{Diff}(S^1)$), the tangent space is isomorphic to the set of all normal $C^\infty$ vector fields along $c$, i.e.,
\begin{equation}
T_\Gamma B_e\cong\{h\colon h=\alpha n,\, \alpha\in C^\infty(S^1,\R)\}
\end{equation}
where $n$ is the unit exterior normal field of the shape $\Gamma$ such that $n (\theta)\perp \Gamma^\prime(\theta)$ for all $\theta\in S^1$ and $\Gamma^\prime$ denotes the circumferential derivative as in \cite{MM-2006}. 
Several intrinsic metrics are discussed in \cite{MM-2006}, among which the following Sobolev metric is used in most algorithmic approaches to shape optimization based on the shape calculus. For $A>0$, the Sobolev metric is induced by the scalar product
\begin{equation}\label{eq_sobolev_metric}
\begin{split}
g^1\colon T_\Gamma B_e\times T_\Gamma B_e & \to \R,\\
    (h,k) & \mapsto
    \int_{\Gamma}\alpha\beta+A\alpha^\prime\beta^\prime
    ds =((id-A\Delta_\Gamma)\alpha,\beta)_{L^2(\Gamma)}
\end{split}
\end{equation}
where $h=\alpha n$ and $k=\beta n$ denote two elements from
the tangent space at $\gamma$ and $\Delta_\Gamma$ denotes the Laplace-Beltrami operator on the surface $\Gamma$.  In \cite{MM-2006} it is shown that the condition $A>0$ guarantees that the scalar product $g^1$ defines a Riemannian metric on $B_e$ and 
thus, geodesics can be used to measure distances.

In \cite{SIOPT2015}, the following scalar product $g^S$ on the tangent space is proposed.
\begin{equation}\label{scp}
\begin{split}
g^S\colon H^{1/2}(\Gamma)\times H^{1/2}(\Gamma) & \to \R,\\
(\alpha,\beta) &\mapsto \langle\alpha,(S^p)^{-1}\beta\rangle=
\int_\Gamma \alpha(s)\cdot [(S^p)^{-1}\beta](s)\ ds.
\end{split}
\end{equation}
where the symmetric and coercive operator $S^p$ is defined by
\begin{equation}
\begin{split}
S^p\colon  H^{-1/2}(\Gamma)& \to H^{1/2}(\Gamma),\\
         \alpha & \mapsto U^\top n
\end{split}
\end{equation}
and $U\in H^1_0(\Omega,\R^d)$ solves the Neumann problem
\begin{equation}\label{weak-elasticity-N2}
a(U,V)=\int_\Gv \alpha\cdot V^\top n\ ds\, , \ \forall\  V\in H^1_0(\Omega,\R^d)
\end{equation}
and $a$ is a coercive and symmetric bilinear form, defined, e.g. by the elasticity equation and thus corresponds to an elliptic problem with fixed outer boundary and forces $\alpha\cdot n $ at the inner boundary $\Gi$.

With the shape space $B_e$ and and a Riemannian metric on its tangent space in hand we can form the Riemannian shape gradient as a Riesz representation of a shape derivative given in the form
\begin{equation}
dJ(\Omega)[V]=\int_\Gamma \gamma\scp{V,n}ds.
\end{equation}
In our model setting the objective function $J$ is given in (\ref{eq_objective}) and its shape derivative in (\ref{obj-der}).
The Riemannian shape gradient $\text{grad}J$ with respect to a Riemannian metric $g\in\{g^1,g^S\}$ is then obtained by
\begin{equation}
\grad J=v n \mbox{ with } g(v,\alpha)=\int_\Gamma\gamma(s)\alpha(s)ds\, , \ \forall \alpha n\in T_\Gamma B_e
\end{equation}
The metric $g^1$, which is also used in \cite{Schulz-Structure-2014} and in many other publications, necessitates a shape derivative in Hadamard form as well as  efficient means to solve linear systems involving the Laplace Beltrami operator in surfaces. All of that requires computational overhead. Furthermore, the surrounding mesh for the computation of the Stokes flow has to be deformed according to the geometry change. This mesh deformation is typically performed by the solution of an elasticity or Poisson problem, with the geometry step as a Dirichlet condition. 

In contrast to that, the usage of the metric $g^S$ corresponds to interpreting the shape derivative as a -- volumetric or boundary -- force to the mesh deformation process as illustrated in \cite{SIOPT2015}. Thus, not only overhead is saved, but also better overall mesh properties are obtained. This paper is devoted to the detailed computational comparison of both metrics. 

Furthermore, it is known that the Riemannian manifold $(B_e, g^1)$ is not metrically complete and the solution of our test problem is not contained in it. On the other hand, according to \cite{SIOPT2015},  the metric $g^S$ gives rise to the manifold 
\begin{equation}\label{shape_manifold}
\SM:=
{\cal H}_i^{1/2}(\Gamma_0,\R^d)\big\slash \mbox{Homeo}(\Gamma_0)
\end{equation}
where $\mbox{Homeo}(\Gamma_0)$ denotes the set of all homeomorphisms of the prior shape $\Gamma_0$ and
\begin{equation}
{\cal H}_i^{1/2}(\Gamma_0,\R^d):=
\{W(\Gamma_0)\colon W\in H^1(\Omega,\Omega),\ W\mbox{ invertible}\}.
\end{equation}
Thus, the construction of $\SM$ is in complete analogy to the construction in \cite{MM-2006} and obviously $B_e\subset \SM$, if $\Gamma_0$ is smooth.
This larger shape manifold (\ref{shape_manifold}) contains the optimal solution of our test problem and is metrically complete.

\section{Augmented Lagrangian method for shape optimization}\label{sec_augmented_lagrange}

The general problem formulation in section 2 includes an objective functional together with a system PDE and geometric constraints. Shape optimization problems treated by means of the shape calculus are not often set within a framework of additional constraints. In our test case, however, the geometric constraints are necessary in order to obtain nontrivial solutions: without the volume constraint $c_{d+1}$, the shape would shrink to a straight line as a trivial but not interesting solution, and without the barycenter constraints $c_{1,\ldots, d}$, the shape would just float out of the computational domain, which again is not a desirable solution. The focus of this paper is on the comparison of the effects of metrics on shapes computed by an optimization process. Thus, we try to keep the optimization framework conceptually as simple as possible and dispense with the potential of additional algorithmic efficiency of one-shot methods \cite{AIAA-2013,BS-CSE-2012,Schulz-LN-2014}, but rather stick with the conceptually simpler and otherwise widely used black-box approach. 
Therefore, we exploit the assumption that the model equation $B_\Omega(y)=0$ can be solved uniquely for $y$, if the shape $\Omega$ is given, which itself depends uniquely on the variable boundary part $\Gamma$. Thus, we consider the system solution $y$ as function of $\Gamma$, i.e., $y=y(\Gamma)$, which results in the following reduced problem formulation:
\begin{align}
\label{red-obj}
\min\limits_\Gamma &J(y(\Gamma))\\
\label{red-constr}
\mbox{s.t.}\ c(\Gamma)&=0
\end{align} 
Indeed, the shape derivative in equation (\ref{obj-der}) is already the shape derivative of (\ref{red-obj}) with respect to $\Gamma$ obeying the chain rule via the implicit function theorem. We note that in this particular self-adjoint combination of PDE and objective we do not need any adjoint equation, which is otherwise usually necessary.

In principle, the optimization problem (\ref{red-obj}, \ref{red-constr}) can be treated by a sequential quadratic programming (SQP) approach on shape manifolds. This would result in the necessity to compute not only the Riesz representation of the shape derivative of the objective, but also of the $(d+1)$ constraints in each optimization step, which increases the algorithmic complexity significantly, since a PDE -- on the surface ($g^1$) or in the volume ($g^S$) -- has to be solved for each of these Riesz representations. This problem can be circumvented by an augmented Lagrangian approach based on the so-called augmented Lagrangian (cf.~\cite{conn1992lancelot}) 
\[
\LL_A(\Gamma,\lambda):=J(y(\Gamma))+\lambda^\top c(\Gamma)+
\frac{\mu}{2}c(\Gamma)^\top c(\Gamma)
\]
for some Lagrange multipliers $\lambda\in\R^{d+1}$ to be determined iteratively below and some $\mu>0$ sufficiently large. The augmented Lagrangian approach has already been used in \cite{EHS2007} in the context of shape calculus. Since the convergence for $\lambda$ is very fast (cf.~Fig.\ref{fig_shape_convergence}), we concentrate in the numerical results section on the local convergence in step 1. for given true multipliers $\lambda$.

Thus, we use the following algorithmic  outline for both metrics to be compared.

\begin{itemize}
\item[0.] Initialize $\lambda^1$, $k:=1$, choose tolerance $\delta_J$ for the optimization process of the augmented Lagrangian and tolerance $\delta_c$ for the satisfaction of constraints $c$. Choose penalty increment factor $\mu_\text{inc}>1$
\item[1.] Solve $\Gamma_k :=\mathop{\mbox{argmin}}\limits_\Gamma \LL_A(\Gamma,\lambda)$ up to tolerance $\delta_J$
\item[2.] If $\Vert c(\Gamma_k)\Vert > \delta_c$
\begin{itemize}
	\item[a)] Update $\mu \leftarrow \mu_\text{inc} \mu$ and go to 1.
\end{itemize}
else
\begin{itemize}
	\item[b)] Update $\lambda^{k+1} \leftarrow \lambda^k+\mu c(\Gamma_k)$ and $k\leftarrow k+1$
\end{itemize}
\item[3.] If $\lambda^k$ is not converged, go to 1.
\end{itemize}

The inner optimization step 1. is chosen as either steepest descent method or as limited memory BFGS-quasi-Newton method as described in \cite{Schulz-Structure-2014} for the metric $g^1$ and in \cite{SIOPT2015} for the metric $g^S$.
We note that the discussion in section \ref{numres} focusses on this step 1.

For the sake of completeness, we rephrase these optimization strategies in the current framework. 
The quasi-Newton approaches rely on the secant condition on manifolds as in \cite{Absil-book-2008}
for a step $\Gamma_{j+1}:=R_{\Gamma_j}(\eta)$ resulting from an increment $\eta_j\in T_{\Gamma_j}\SM$ in iteration $j$ via a retraction $R$ as
\[
\grad \LL_A (\Gamma_{j+1})-{\cal T}_{\eta_{j}}\grad \LL_A (\Gamma_{j})=G_{j+1}[{\cal T}_{\eta_{j}}\eta_{j}]
\]
where ${\cal T}:T\SM\oplus T\SM\to T\SM: (h_\Gamma,k_\Gamma)\mapsto {\cal T}_{h_\Gamma}k_\Gamma$ is a vector transport associated to the retraction $R$ and $G_{j+1}$ is intended to approximate the
Riemannian Hessian $\nabla\grad \LL_A(\Gamma_{j+1})$.

In the following limited BFGS loop, we use the notation
\begin{align*}
s_j:=&{\cal T}_{\eta_{j}}\eta_{j}\in T_{\Gamma_{j+1}}\SM\\
y_j:=& \grad \LL_A (\Gamma_{j+1})-{\cal T}_{\eta_{j}}\grad \LL_A (\Gamma_{j})\in T_{\Gamma_{j+1}}\SM
\end{align*}
In \cite{Ring-Wirth-2012}, superlinear convergence properties for BFGS-quasi-Newton methods on manifolds are analyzed for the case that ${\cal T}_{\eta_j}$ is
an isometry. This requirement is satisfied, e.g., if $\cal T$ and $R$
are the parallel transport and the exponential map. Details on the specific operators $\mathcal{T},R$ used here are given in section \ref{numres}.

Thus steepest descent methods and quasi-Newton methods, used here, are described jointly by the following algorithm

\medskip

\begin{algorithm}[H]
 \caption*{\bf Details of step 1.}
 \label{algopt}
\begin{algorithmic}
\Repeat
\State $\rho_j \gets g(y_j,s_j)^{-1}$
\State $q \gets \grad \LL_A(\Gamma_j)$
\For{$i = j-1, \dots , j-m$}
	\State $s_i \gets {\cal T}_{q}s_i$
	\State $y_i \gets {\cal T}_{q}y_i$
	\State $\alpha_i \gets \rho_i g(s_i,q)$
	\State $q \gets q - \alpha_i y_i$
\EndFor
\State $q \gets \frac{g(y_{j-1},s_{j-1})}{g(y_{j-1},y_{j-1})} q$
\For{$i = j-m, \dots , j-1$}
	\State $\beta_i \gets \rho_i g(y_i,q)$
	\State $q \gets q + (\alpha_i - \beta_i) s_i$
\EndFor 
\State $\Gamma_{j+1} \gets R_{\Gamma_j}(q)$\\
\Until $\Vert q \Vert_2 < \delta_J$
\end{algorithmic}
\end{algorithm}

If $m=0$ above, the algorithm boils down to steepest descent methods. Note that the algorithm only describes the deformation of the shape $\Gamma$. The surrounding mesh is deformed according to an elasticity equations in the case $g=g^1$ and in the case $g=g^S$ by the usage of the elastic deformation field which is anyway available from the computation of $q$, such that the shape derivative is interpreted as a boundary force in the latter case.

\section{Numerical results}\label{numres}

\begin{figure}
\begin{center}
\includegraphics[width=0.6\textwidth]{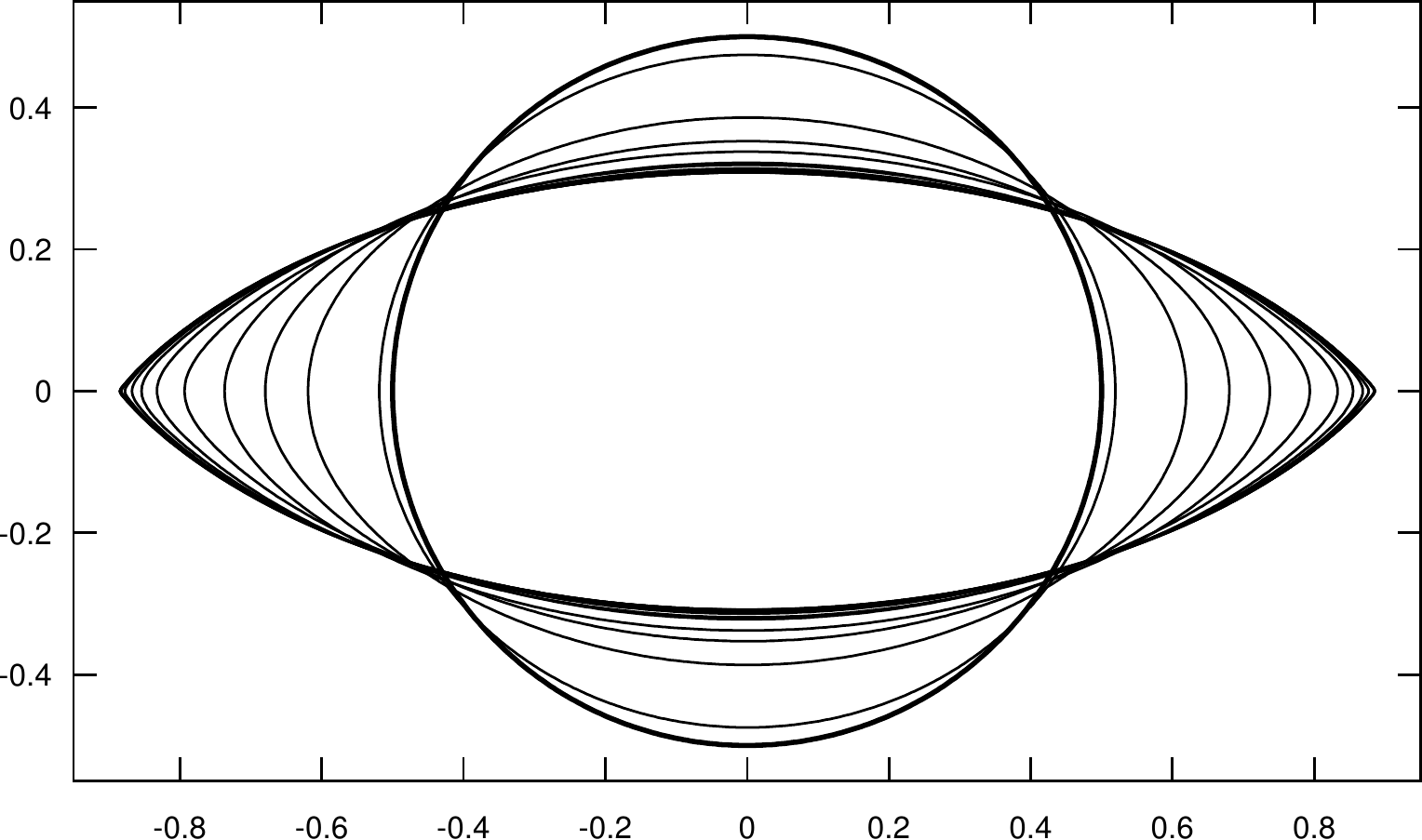}
\end{center}
\caption{Iterations of the BFGS method with Steklov-Poincaré metric (initial and optimal shapes are bold)}
\label{fig_bfgs_shapes}
\end{figure}

We are now prepared to describe a specific test case for the algorithms outlined in the previous sections.
The computational domain, as depicted in figure \ref{fig_schematic}, is chosen to be $\Omega_0 \cup \Oe = [-3, 6]\times[-2,2]$ for the 2 dimensional case.
The initial body $\Omega_0$ is a circle with barycenter $\text{bc}(\Omega_0) = (0,0)$ and radius $r= 0.5$ leading to $\text{vol}(\Omega_0) = \frac{\pi}{4}$.
The computational grid consists of 10,150 triangles of which 633 form the variable surface $\Gv$.
Analogously, the 3 dimensional domain is generated by rotating the 2d mesh around the X axis.
Here we have 27,892 tetrahedral elements with 1,206 triangles forming $\Gv$.

\begin{figure}
\begin{center}
\begin{subfigure}{0.49\textwidth}
\includegraphics[width=1.0\textwidth]{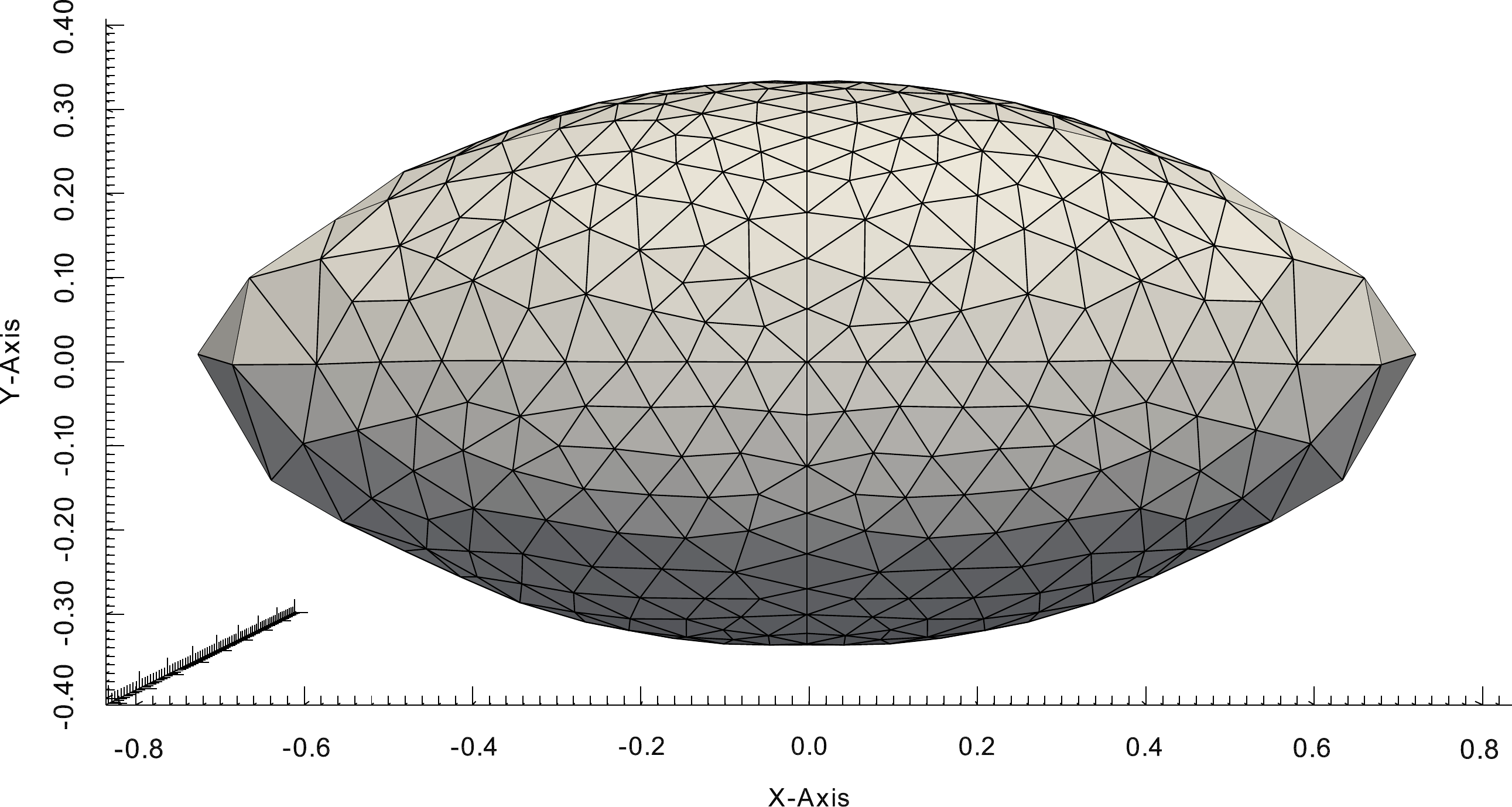}
\subcaption{Laplace-Beltrami metric, no FEM solution possible}
\label{fig_shape_3d_LB}
\end{subfigure}
\begin{subfigure}{0.49\textwidth}
\includegraphics[width=1.0\textwidth]{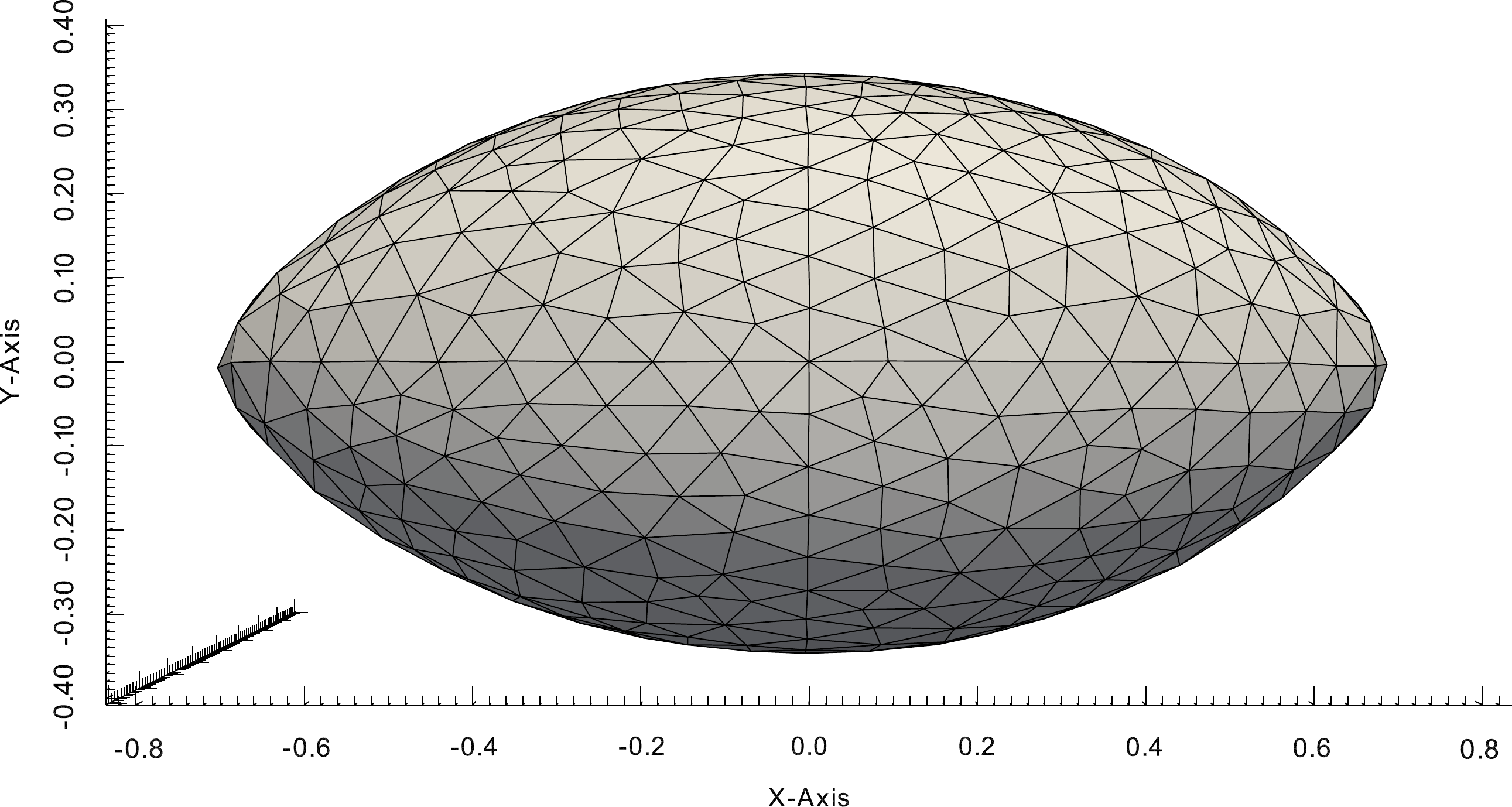}
\subcaption{Steklov-Poincaré metric, converged solution}
\label{fig_shape_3d_SP}
\end{subfigure}
\end{center}
\caption{3D comparison if $g^1$ after 62 gradient steps and $g^S$ metric with BFGS updates after 16 iterations}
\label{fig_shape_3d}
\end{figure}

Compared to the 3d mesh we have chosen a much finer discretization for the 2d case.
This is due to the fact that we want to measure convergence of the proposed algorithm which is only practical in 2d.
Whereas, the coarse 3d discretization shall demonstrate the ability of the shape metric $g^S$ to also handle sharp edges in this situation.
This can be seen in figure \ref{fig_shape_3d}.
Here we again observe that the $g^S$ metric is superior compared to $g^1$ with respect to the node distribution on the surface of the body.
While the solution on the right hand side is converged to the actual solution, the $g^1$-algorithm on the left hand side breaks down too early with unfeasible grids.

It should be remarked that there is no mesh inside $\Omega$.
Thus, we apply divergence theorem to the constraints $c$ and obtain
\begin{equation}
\text{vol}(\Omega) = -\int_\Gv \frac{1}{d} \langle s , n \rangle\, ds \quad \text{and} \quad \text{bc}(\Omega) = -\frac{1}{2\text{vol}(\Omega)}\int_\Gv \langle (s_1^2, \dots, s_d^2)^T, n \rangle\, ds.
\end{equation}
One could also constrain the volume and barycenter of $\Oe$ to be constant which is equivalent since the outer boundaries are fixed.

From a computational point of view, the algorithms for the $g^S$ metric are favorable compared to the $g^1$ metric.
In any case, we discretize the velocity $v$ with continuous piece-wise quadratic (P2) elements and the pressure $p$ with continuous piece-wise linear (P1) elements in order to guarantee stability, which is the standard approach for Stokes equation.
This discretization leads to a major advantage of the $g^S$ metric.
Since $v$ is represented in P2 functions, we end up with a discontinuous piece-wise linear representation for the term $\sum_{i=1}^d \left( \frac{\partial v_i}{\partial n} \right)^2$ in equation \eqref{obj-der}.
However, parts of the shape derivative coming from the geometric constraints are piece-wise linear due to the P1 shape functions of the straight-lined elements.
In the $g^1$ case we thus have to perform $L^2$-projections in order to represent all quantities in the same basis functions.
Given a function $u \in V_1$ in a finite element space which has to be projected into a different space $V_2$ we have to solve
\begin{equation}
\int_\Gv \bar{u} v\, ds = \int_\Gv u v\, ds\quad \forall\, v \in V_2
\end{equation}
for $\bar{u}$.
Whereas in the $g^S$ case these terms only show up on the right hand side of equation \eqref{weak-elasticity-N2}.
It is thus straightforward to combine expressions which are represented in different basis functions.
By solving the Neumann problem \eqref{weak-elasticity-N2}, we obtain a representation of the shape derivative with respect to the $g^S$ metric.
However, the solution $U$ of \eqref{weak-elasticity-N2} is not only defined on $\Gv$ but also in the entire domain $\Oe$.
This gives us a smooth deformation field which can be applied to the finite element mesh.
In contrast, using the $g^1$ metric we have to solve a tangential Laplace equation yielding a representation of the shape gradient only defined at $\Gv$.
In an additional step this information is plugged into a linear elasticity problem as Dirichlet boundary condition.

Comparing the computational costs of the two algorithms we additionally have to solve two PDEs defined at $\Gv$ in the $g^1$ case.
This stems from the fact that gradient representation and mesh deformation are computed in one step by using the $g^S$ metric.

\begin{figure}
\begin{center}
\begin{subfigure}{0.49\textwidth}
\includegraphics[width=1.0\textwidth]{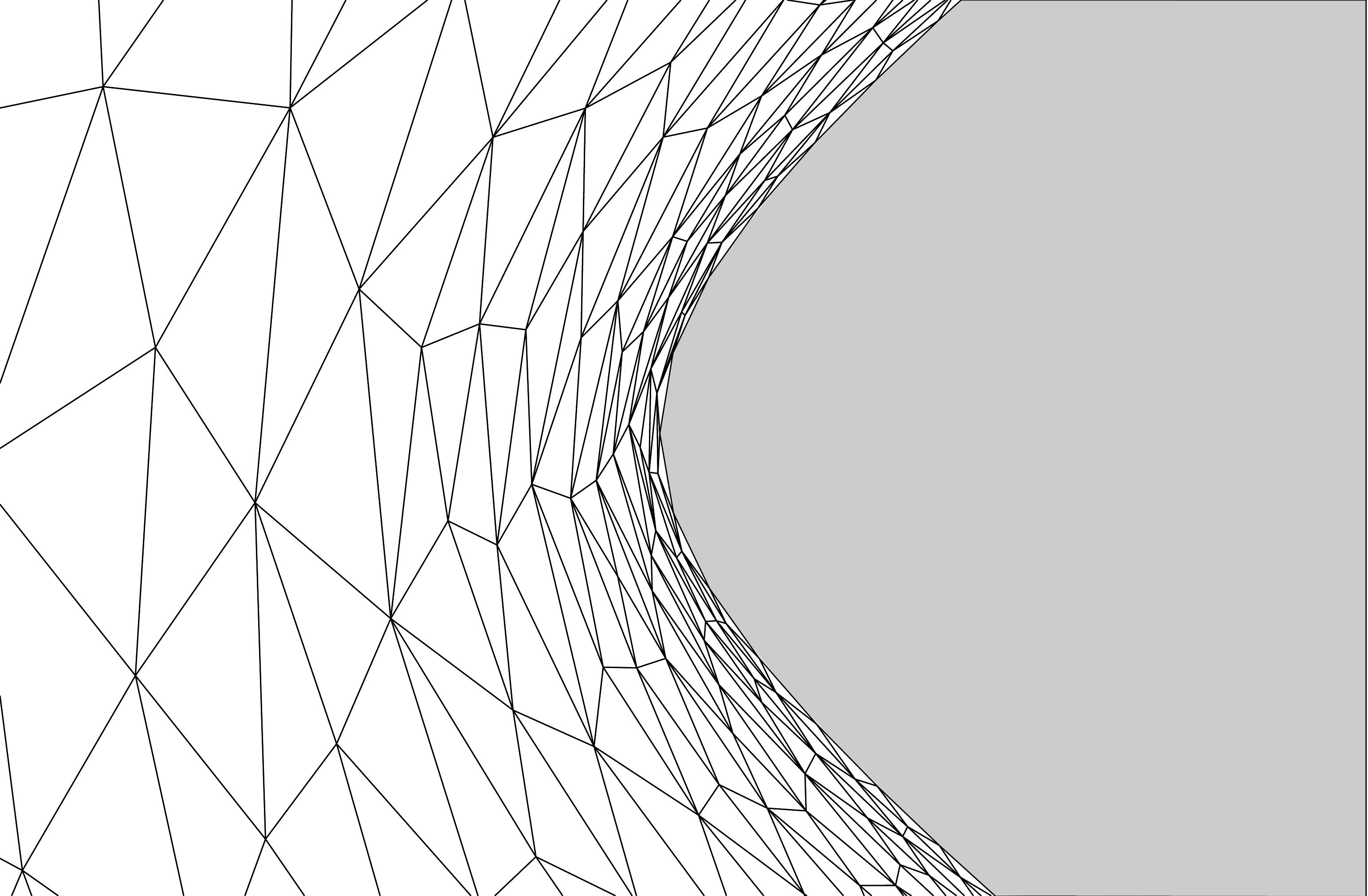}
\subcaption{Laplace-Beltrami metric, no FEM solution possible}
\label{fig_optimal_mesh_LB}
\end{subfigure}
\begin{subfigure}{0.49\textwidth}
\includegraphics[width=1.0\textwidth]{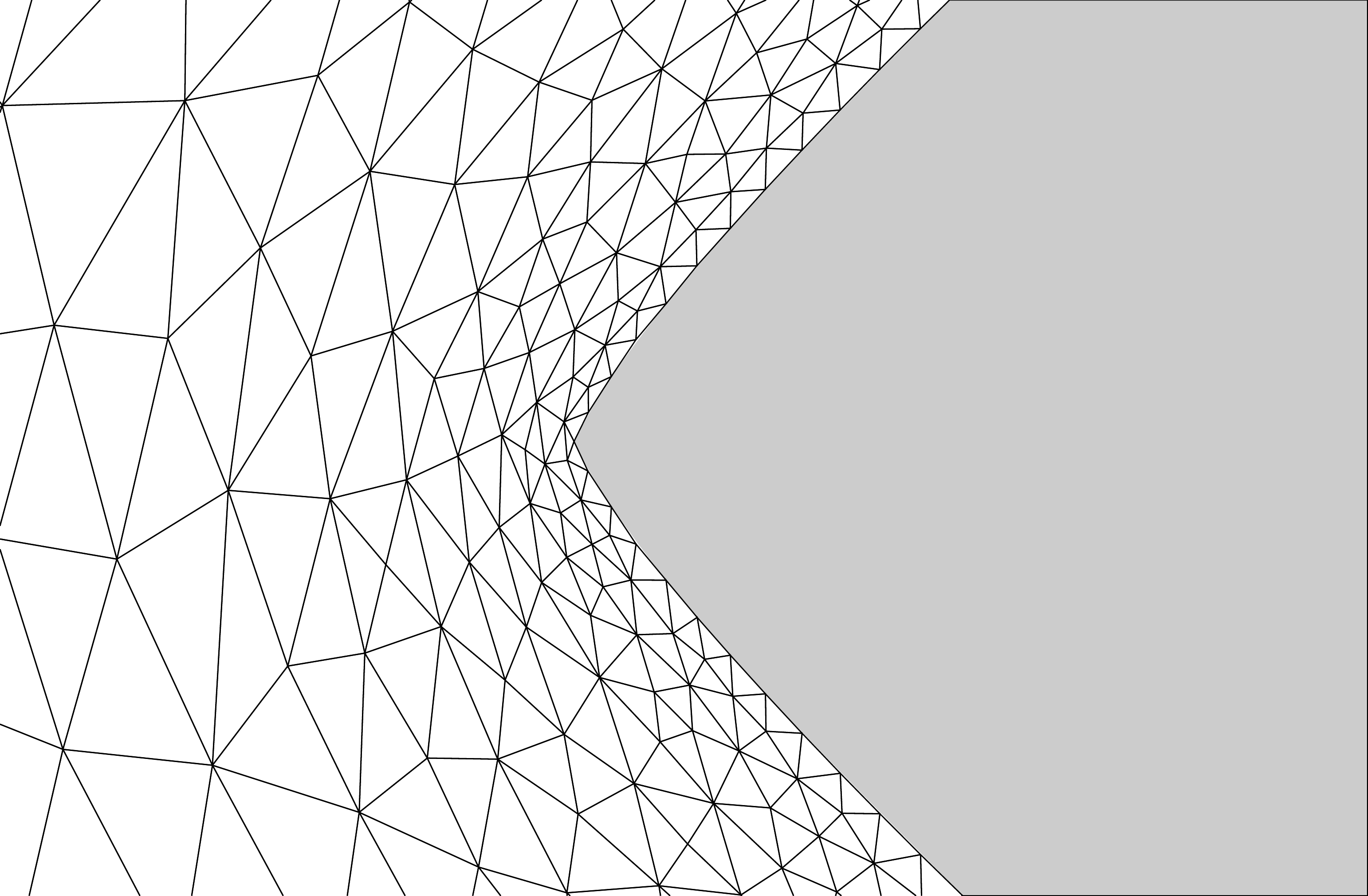}
\subcaption{Steklov-Poincaré metric, converged solution}
\label{fig_optimal_mesh_SP}
\end{subfigure}
\end{center}
\caption{Visual comparison of shape metrics with respect to influence on mesh quality (tip of the body), left: 65 gradient steps, right: 16 BFGS steps}
\label{fig_metrics_comparison}
\end{figure}

Secondly, what makes the Steklov-Poincaré metrics an even more appealing concept is the mesh quality that results from this two approaches.
The classical approach, which is based on the $g^1$ metric, only takes deformations normal to $\Gv$ into account.
This reflects the Hadamard theorem stating that only the normal component of deformations affect the objective function.
However, using metrics of Steklov-Poincaré type also allows surface nodes to slide along $\Gv$ which strongly influences the mesh quality.
This can be seen in figure \ref{fig_metrics_comparison}.
In this particular situation, a solution of the optimization problem \eqref{eq_objective} subject to \eqref{eq_stokes} and \eqref{eq_geo_constraints} can not be achieved with the classical approach.
These iterations lead to discretization meshes that do not allow further computations due to degenerated cells.

\begin{figure}
\begin{center}
\includegraphics[width=0.8\textwidth]{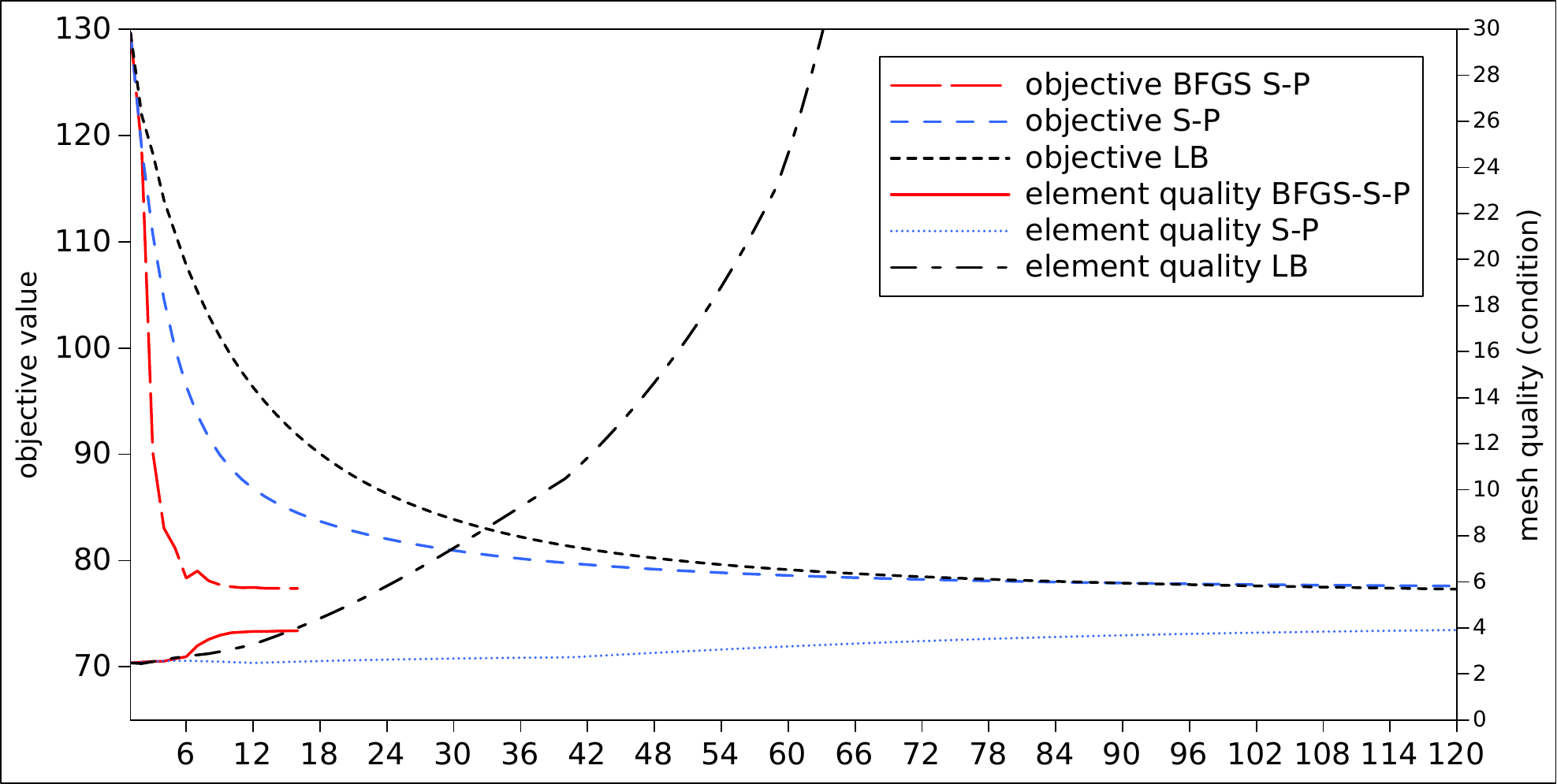}
\end{center}
\caption{Objective values and mesh quality for Steklov-Poincaré and Laplace-Beltrami metric}
\label{fig_mesh_quality}
\end{figure}

Figure \ref{fig_mesh_quality} shows the decreasing mesh quality during the optimization iterations.
Here we assume the Lagrange multipliers to be known and start the optimization with the circle geometry.
The mesh quality is measured with respect to the condition of the affine mapping between reference and physical element.
Values close to one indicate a good mesh quality.
Here the worst element is shown.
For the Steklov-Poincaré type metric we achieve mesh qualities which are within $[2.47, 3.87]$ during the entire optimization.
Whereas, the classical approach leads to degenerated elements, which can be observed by the unbound condition numbers.
This behavior makes it impossible to apply BFGS updates to the Laplace-Beltrami-based algorithm as the algorithm breaks down in the first iterations with unfeasible elements due to the larger step sizes.
In contrast, figure \ref{fig_mesh_quality} also shows the speedup that can be gained by with a quasi Newton method.
In this particular case we applied a limited memory BFGS strategy with 3 gradients in storage.
A more detailed description of quasi Newton updates for shape optimization can be found in \cite{SIOPT2015}.
The corresponding 16 iterated shapes are shown in figure \ref{fig_bfgs_shapes}.

The key part of the method described here is the choice of the bilinear form in \eqref{weak-elasticity-N2}.
It turns out that the mesh quality heavily depends on $a$.
The results shown in this work are obtained by choosing $a$ as the weak form of the linear elasticity equations
\begin{equation}\label{eq_linelas}
\begin{aligned}
\text{div}( \sigma ) &= 0 &\quad& \text{in} \quad \Oe\\
U &= 0 &\quad& \text{on} \quad \Go \cup \Gi \cup \Gw\\
\frac{\partial U}{\partial n} &= \gamma &\quad& \text{on} \quad \Gv\\
\end{aligned}
\end{equation}
in terms of 
\begin{equation}
\begin{aligned}
\sigma &:= \lambda_\text{elas} \text{Tr}(\epsilon) I + 2 \mu_\text{elas} \epsilon\\
\epsilon &:= \frac{1}{2}(\nabla U + \nabla U^T)\\
\end{aligned}
\end{equation}
where $\sigma$ and $\epsilon$ are the strain and stress tensor, respectively.
Here $\lambda_\text{elas}$ and $\mu_\text{elas}$ denote the Lam\'{e} parameters, which can be expressed in terms of Young's modulus $E$ and Poisson's ratio $\nu$ as
\begin{equation}
\lambda_\text{elas} = \frac{\nu E}{(1+\nu)(1-2\nu)} \, , \quad \mu_\text{elas} = \frac{E}{2(1+\nu)}.
\end{equation}
Equation \eqref{weak-elasticity-N2} then transforms to
\begin{equation}
\int_\Oe \sigma(U):\epsilon(V) \, dx = \int_\Gv \gamma \langle n, V \rangle\, ds \quad \forall \;V \in H_0^1 (\Omega, \mathbbm{R}^d).
\end{equation}
The optimal mesh shown in figure \ref{fig_optimal_mesh_SP} is obtained with $\lambda_\text{elas} = 0$ and $\mu_\text{elas} \in [\mu_\text{min}, \mu_\text{max}] = [1, 500]$ smoothly decreasing from $\Gv$ to the outer boundaries.
We therefor solve Poissons equation
\begin{equation}
\begin{aligned}
\Delta \mu_\text{elas} &= 0 &\quad& \text{in} \quad \Oe\\
\mu_\text{elas} &= \mu_\text{max} &\quad& \text{on} \quad \Gv\\
\mu_\text{elas} &= \mu_\text{min} &\quad& \text{on} \quad \Go \cup \Gi \cup \Gw
\end{aligned}
\end{equation}
in an initial stage of the optimization.
It should be remarked that the solution behavior of Poissons equation depends on the dimension.
In the 3d case we therefor apply $\sqrt{\mu_\text{elas}}$ in contrast to $\mu_\text{elas}$ for the elasticity tensor of the $g^S$ metric.

Still two questions remain open, namely the numerical realization of the operators $\mathcal{T}$ and $R$.
In principle, the retraction $R$ would require yet another solution of a PDE, which is computationally too expensive.
As discussed in \cite{Schulz-Structure-2014}, simply adding the computed deformation field $U$ to the nodes of the finite element mesh approximates the retraction in a reasonable way.
Similarly, the vector transport necessary for the limited BFGS algorithm is approximated with the identity operator.
From a computational point of view, $s_i$ and $y_i$ are finite dimenional vectors in the memory.
Thus, we use these vectors as they are, without taking into account that they approximate elements in different tangent spaces.

\begin{figure}
\begin{center}
\includegraphics[width=0.4\textwidth]{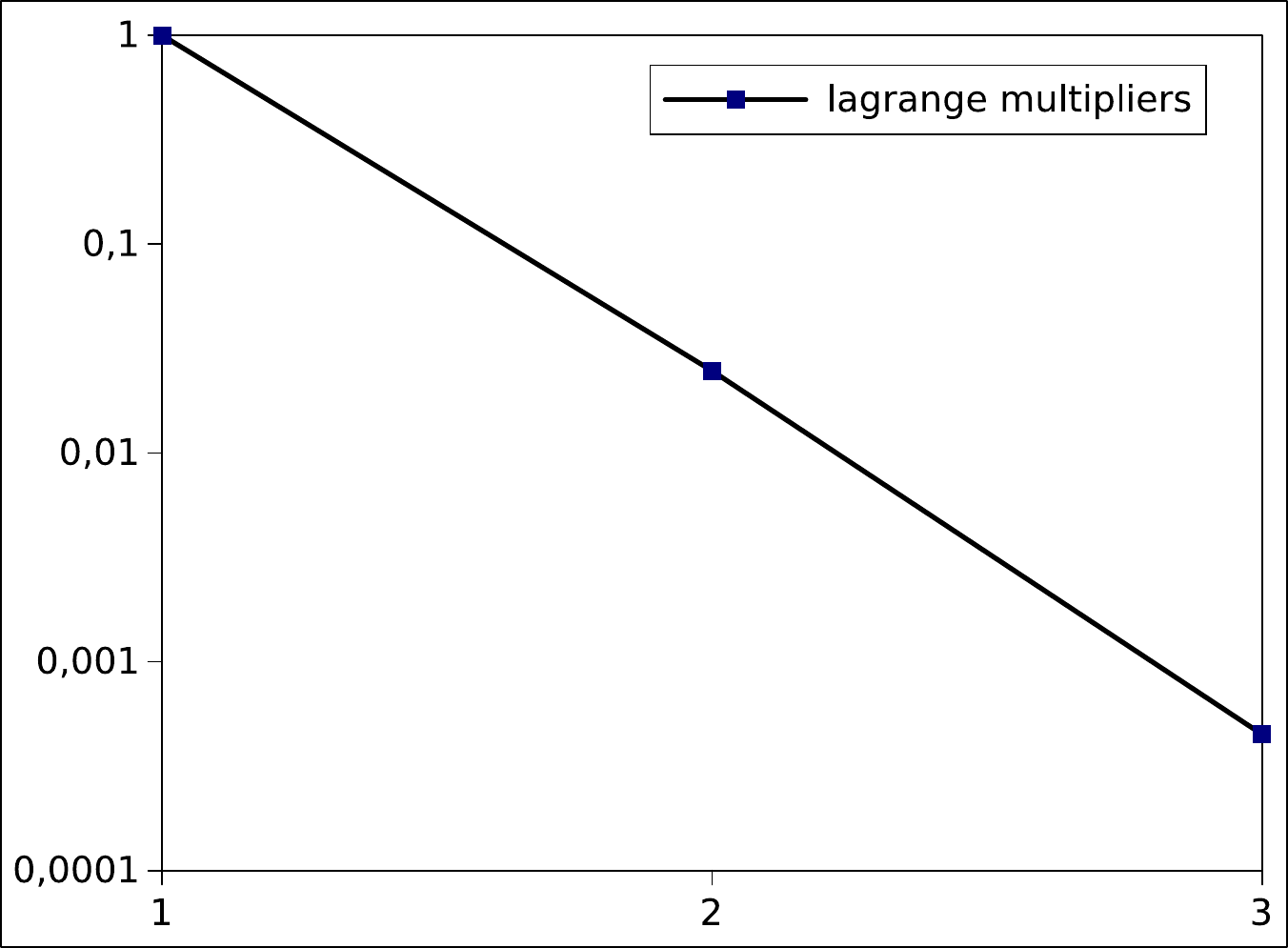}
\end{center}
\caption{Convergence of the Lagrange multipliers for the geometric constraints. The figure shows the distance to the final Lagrange multiplier in $L^2$-norm. }
\label{fig_lambda_convergence}
\end{figure}

\begin{figure}
\begin{center}
\includegraphics[width=0.4\textwidth]{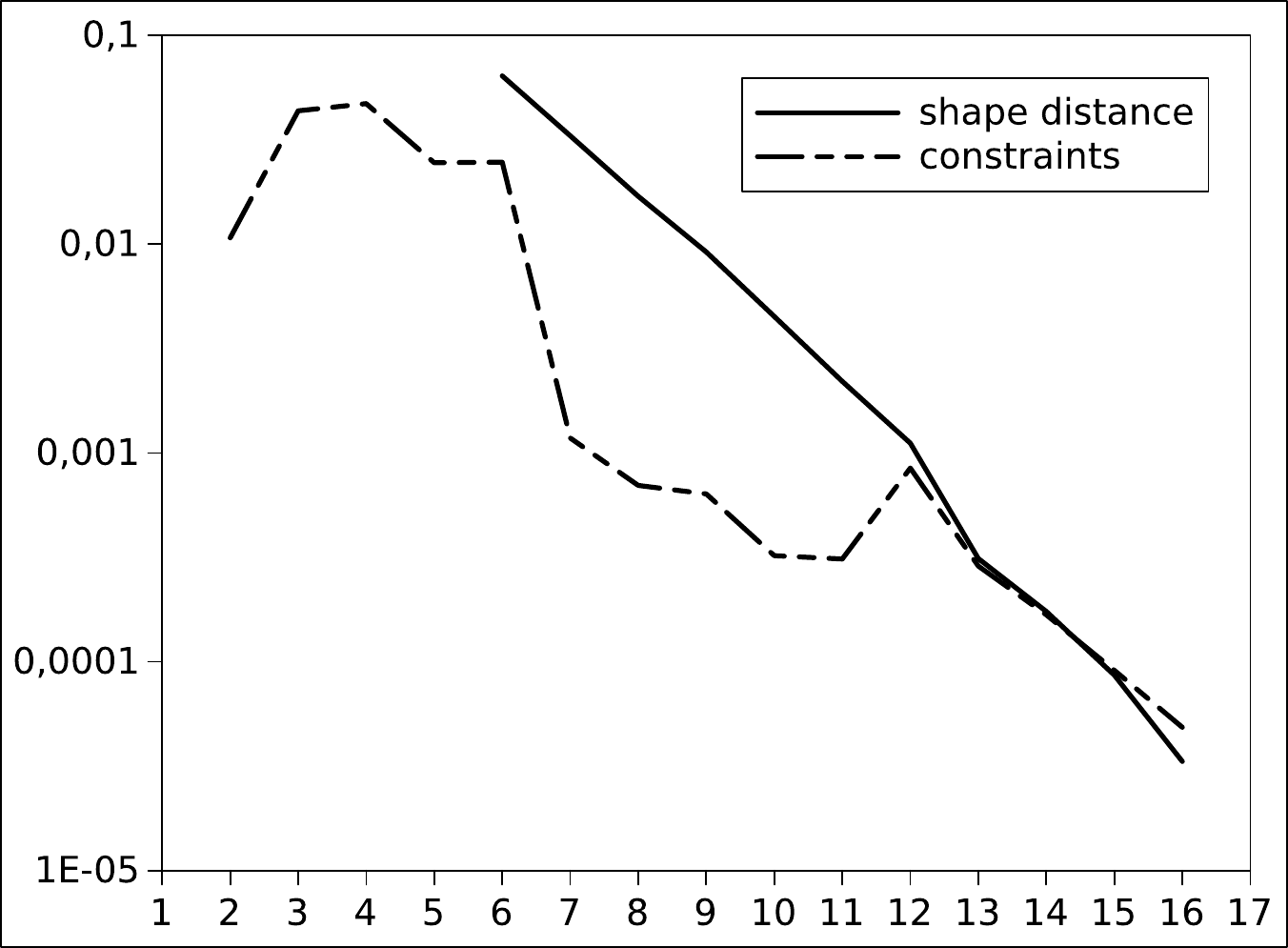}
\end{center}
\caption{Shape distances to optimal solution and convergence of geometric constraints for BFGS-Steklov-Poincaré, both with fixed Lagrange multipliers}
\label{fig_shape_convergence}
\end{figure}

In order to estimate the convergence speed of the proposed method we have to use a measure for the distance between two shapes.
We therefore approximate the geodesic distance of the current to the converged shape.
This is done by finding the pointwise distance in normal direction between the two shapes and integrating this quantity over $\Gv$.
This can be seen in figure \ref{fig_shape_convergence}.
Here the convergence of the corresponding constraints $\Vert c(\Gv^k) \Vert_2$ is also visualized.
The underlying optimization is performed with the converged Lagrange multipliers for the geometric constraints.
Note that values of the shape distance for the first iterations are missing.
This is due to the fact that the underlying shape metric as described above is only reasonable for small deformations and thus does not work for the first iterations.

As mentioned in section \ref{sec_augmented_lagrange}, the convergence of $\lambda$ is relatively fast.
This can be seen by observing $\Vert \lambda \Vert_2$ which is shown in figure \ref{fig_lambda_convergence}.
For this test we have chosen the tolerance for the inner optimization loop depending on the $L^2$ norm of the step $U$ evaluated on $\Gv$ to be $\delta_J = 10^{-4}$.
These computations are performed with a penalty factor $\mu = 10^2$.
In principle, the augmented Lagrangian algorithm provides updates for the penalty factor, if the violation of the constraints is over a given tolerance $\delta_c$.
We yet did not encounter convergence problems in the test cases with the penalty factor $\mu = 10^2$, since, as figure \ref{fig_shape_convergence} suggests, a convergence of the constraints is achieved after appropriate Lagrange multipliers are found by the algorithm.

The algorithms presented within this work are entirely implemented using the GetFEM++ library for the assembly of the finite elements and the PETSc library for the solution of the linear systems.
The initial grids are generated with the Gmsh mesh generator and the mesh optimization routines therein.

\section{Conclusions}
We compare computational aspects of the standard surface metric for shape optimization based on shape calculus with the surface metric introduced and analyzed in \cite{SIOPT2015}. Advantages of the latter metric are demonstrated with regards to convergence properties, computational overall effort as well as resulting mesh quality. This holds for a two-dimensional and for a three-dimensional set-up.

\section*{Acknowledgment}
This work has been partly supported by the Deutsche Forschungsgemeinschaft within the Priority program SPP 1648 ``Software for Exascale Computing'' under contract number Schu804/12-1. Furthermore, we acknowledge support by the  Sino-German Science Center  on the occasion of the Chinese-German Workshop on Computational and Applied Mathematics in Augsburg 2015.

\bibliographystyle{plain}
\bibliography{citations.bib}

\begin{thebibliography}{10}

\bibitem{Absil-book-2008}
{P.-A.} Absil, R.~Mahony, and R.~Sepulchre.
\newblock {\em Optimization Algorithms on Matrix Manifolds}.
\newblock Princeton University Press, 2008.

\bibitem{optbook}
A.~Borz\`{\i} and V.~H. Schulz.
\newblock {\em Computational optimization of systems governed by partial
  differential equations}.
\newblock Number~08 in SIAM book series on Computational Science and
  Engineering. SIAM Philadelphia, 2012.

\bibitem{BS-CSE-2012}
Alfio Borz\`{\i} and Volker Schulz.
\newblock {\em Computational Optimization of Systems Governed by Partial
  Differential Equations}.
\newblock Society for Industrial and Applied Mathematics, 2012.

\bibitem{conn1992lancelot}
A.~R. Conn, N.~I.~M. Gould, and Ph.~L. Toint.
\newblock {\em Lancelot}.
\newblock Springer Berlin Heidelberg, 1992.

\bibitem{Delfour-Zolesio-2001}
M.~C. Delfour and J.-P. Zol\'esio.
\newblock {\em Shapes and Geometries: Analysis, Differential Calculus, and
  Optimization}.
\newblock Advances in Design and Control. SIAM Philadelphia, 2001.

\bibitem{EHS2007}
K.~Eppler, H.~Harbrecht, and R.~Schneider.
\newblock On convergence in elliptic shape optimization.
\newblock {\em SIAM J. Control Optim.}, 46(1):61--83, 2007.

\bibitem{Langer-2015}
P.~Gangl, A.~Laurain, H.~Meftahi, and K.~Sturm.
\newblock Shape optimization of an electric motor subject to nonlinear
  magnetostatics.
\newblock Technical report, \url{http://arxiv.org/abs/1501.04752}, 2015.

\bibitem{Haack1941}
W.~Haack.
\newblock Gescho\ss{}formen kleinsten {W}ellenwiderstandes.
\newblock {\em Bericht der Lilienthal-Gesellschaft}, 136(1):14–28, 1941.

\bibitem{MM-2006}
{P. W.} Michor and D.~Mumford.
\newblock Riemannian geometries on spaces of plane curves.
\newblock {\em J. Eur. Math. Soc. (JEMS)}, 8:1--48, 2006.

\bibitem{Mohammadi-2001}
B.~Mohammadi and O.~Pironneau.
\newblock {\em Applied Shape Optimization for Fluids}.
\newblock Numerical Mathematics and Scientific Computation. Clarendon Press
  Oxford, 2001.

\bibitem{Skin-2015}
A.~N\"agel, V.~Schulz, M.~Siebenborn, and G.~Wittum.
\newblock Scalable shape optimization methods for structured inverse modeling
  in 3{D} diffusive processes.
\newblock {\em Computing and Visualization in Science}, 2015.

\bibitem{Paganini}
A.~Paganini.
\newblock Approximative shape gradients for interface problems.
\newblock Technical Report 2014-12, Seminar for Applied Mathematics, ETH
  Z{\"u}rich, 2014.

\bibitem{Pironneau1973}
O.~Pironneau.
\newblock On optimum profiles in stokes flow.
\newblock {\em Journal of Fluid Mechanics}, 59(1):117–128, 1973.

\bibitem{Ring-Wirth-2012}
W.~Ring and B.~Wirth.
\newblock Optimization methods on {R}iemannian manifolds and their application
  to shape space.
\newblock {\em {SIAM} Journal of Optimization}, 22:596--627, 2012.

\bibitem{AIAA-2013}
S.~Schmidt, C.~Ilic, V.~Schulz, and N.~Gauger.
\newblock Three dimensional large scale aerodynamic shape optimization based on
  the shape calculus.
\newblock {\em AIAA Journal}, 51(11):2615--2627, 2013.

\bibitem{Schulz-Structure-2014}
V.~Schulz, M.~Siebenborn, and K.~Welker.
\newblock Structured inverse modeling in parabolic diffusion problems.
\newblock {\em SIAM Control}, 2014.
\newblock \url{http://arxiv.org/abs/1409.3464} (accepted).

\bibitem{Schulz-LN-2014}
V.~Schulz, M.~Siebenborn, and K.~Welker.
\newblock Towards a {L}agrange-{N}ewton approach for {P}{D}{E} constrained
  shape optimization.
\newblock In {\em Trends in PDE Constrained Optimization}, volume 165 of {\em
  International Series of Numerical Mathematics}. Springer, 2014.
\newblock \url{http://arxiv.org/abs/1405.3266}.

\bibitem{SIOPT2015}
V.~Schulz, M.~Siebenborn, and K.~Welker.
\newblock A novel {S}teklov-{P}oinca\'e type metric for efficient {P}{D}{E}
  constrained optimization in shape spaces.
\newblock {\em SIAM Optimization}, 2015.
\newblock \url{http://arxiv.org/abs/1506.02244} (submitted).

\bibitem{VHS-shape-Riemann}
{V. H.} Schulz.
\newblock A {R}iemannian view on shape optimization.
\newblock {\em Foundations of Computational Mathematics}, 14:483--501, 2014.

\bibitem{SokoZol}
J.~Sokolowski and {J.-P.} Zol\'{e}sio.
\newblock {\em An introduction to shape optimization}.
\newblock Springer, 1992.

\bibitem{Berggren-horn-2007}
R.~Udawalpola and M.~Berggren.
\newblock Optimization of an acoustic horn with respect to efficiency and
  directivity.
\newblock {\em Internat. J. Numer. Methods Engrg.}, 73(11):1571--1606, 2007.

\end{thebibliography}

\end{document}